\documentclass[a4,11pt]{amsart}
\usepackage{amssymb,amsthm, amsmath}
\usepackage{graphicx}
\usepackage{subfigure}
\def\N{\mathbb N}
\def\Z{\mathbb Z}

\def\R{\mathbb R}
\def\C{\mathbb C}

\def\T{\mathbb T}

\newtheorem{theorem}{Theorem}
\newtheorem{thm}[theorem]{Theorem}

\newtheorem{lem}[theorem]{Lemma}
\newtheorem{cor}[theorem]{Corollary}

\newcommand{\e}{\mathop{\mathbf{e}}\nolimits}

\begin{document}
\title{Spiral Delone sets and three distance theorem}
\date{}
\author{Shigeki Akiyama}
\begin{abstract}
We show that a constant angle progression on the Fermat spiral forms a 
Delone set if and only if its angle is badly approximable. 
\end{abstract}
\maketitle

In botany, phyllotaxis stands for angular patterns of plant stems or leaves,  
which tends to be correlated with Fibonacci numbers and the golden mean. 
In this short note,  we study a related fundamental 
problem on the spiral patterns, which is not addressed before.
Let $X$ be a subset of $\R^2$ which is identified with the complex plane $\C$.
Denote by $B(x,r)$ the open ball of radius $r$ centered at $x$.
$X$ is $r$-{\it relatively dense} if there exists $r>0$ that for any $x\in \C$,  
$B(x,r)\cap X\neq \emptyset$ holds. 
$X$ is $s$-{\it uniformly discrete}
if there exists $s>0$ that for any $x\in \C$
we have $\text{Card}(B(x,s)\cap X)\le 1$. 
$X$ is a $(r,s)$-{\it Delone set} if it is both $r$-relatively dense and 
$s$-uniformly discrete.
It is clear that $r$-relatively dense implies $r'$-relatively dense if $r\le r'$, and 
$s$-uniformly discrete implies $s'$-uniformly discrete if $s\ge s'$, and $r\ge s$ holds 
for a $(r,s)$-Delone set. 
We omit the prefixes $r$-, $s$- and $(r,s)$- if we are
only interested in the existence of $r$ and/or $s$.

Denote by $\e(z)=e^{2\pi z \sqrt{-1}}$. Fix an {\it angle} $\alpha\in [0,1)$ and a 
strictly increasing function $f(t)$ from $\R_{\ge 0}$ to itself.  We study a point set
$$
X(f,\alpha)=\{ f(n)\e(n \alpha) \ |\  n\in \N \}
$$
on a spiral curve $\{ f(t)\e(t \alpha)\ |\ t\in \R_{\ge 0} \}$.
Clearly $X(f,\alpha)$ is not relatively dense if the angle $\alpha$ is rational, since
$X(f,\alpha)$ is contained in a union of finite number of lines passing through the origin.

\begin{lem}
\label{Ness}
If $X(f,\alpha)$ is $r$-relatively dense, then $\limsup_{n\rightarrow \infty} f(n)/\sqrt{n}<2r$.
If $X(f,\alpha)$ is $s$-uniformly discrete, then $\liminf_{n\rightarrow \infty} f(n)/\sqrt{n}>s/2$. 
\end{lem}

\proof
Assume that $X(f,\alpha)$ is $r$-relatively dense. Then 
$$\bigcup_{z\in X(f,\alpha)} B(x,r)\supset \C.$$ 
Since $K[n]:=\{ z\in X(f,\alpha)\ |\ |z|\le f(n)\}$
has cardinality $n$ and 
$$
\bigcup_{z\in K[n]} B(z,r) \supset B(0,f(n)-r)
$$
for $f(n)>r$, we see $n\pi r^2\ge \pi(f(n)-r)^2$ which implies $2r\ge f(n)/\sqrt{n}$.

If $X(f,\alpha)$ is $s$-uniformly discrete, then $B(z,s)$ are disjoint disks for $z\in X(f,\alpha)$. 
From $K[n]\subset B(0,f(n))$, we obtain 
$$
\bigcup_{z\in K[n]} B(z,s) \subset B(0,f(n)+s)
$$
which leads to
$n \pi s^2 \le \pi (f(n)+s)^2$, i.e., $s/2 \le f(n)/\sqrt{n}$ for $n\ge 4$.
\qed
\bigskip

Our target is to obtain a condition that $X(f,\alpha)$ is a Delone set. 
It is of interest to discuss general increasing functions $f$, however, 
in this paper 
we specify $f(t)=\sqrt{t}$ and $\alpha$ is irrational and study the set
$$
X(\sqrt{t},\alpha)=\{ \sqrt{n} \e(n \alpha) \ |\ n\in \N \}.
$$
From now on, we assume 
\begin{equation}
\label{rs}
r>1/2 \quad \text{ and } \quad s<2
\end{equation}
in light of Lemma \ref{Ness}.

\begin{figure}[htb]
\begin{subfigure}[$X(\sqrt{t},\frac{3-\sqrt{5}}2)$]{ 
\includegraphics[clip,width=0.45\columnwidth]{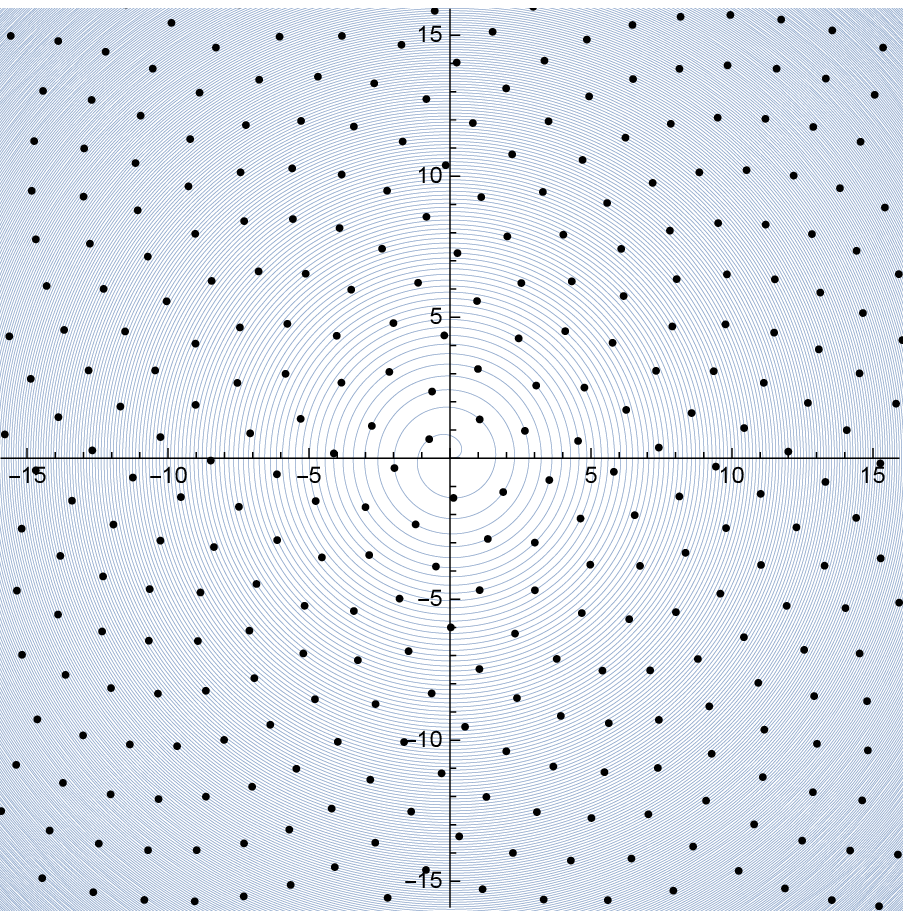}}
\end{subfigure}~
\begin{subfigure}[$X(\sqrt{t},\pi-3)$]{
\includegraphics[clip,width=0.45\columnwidth]{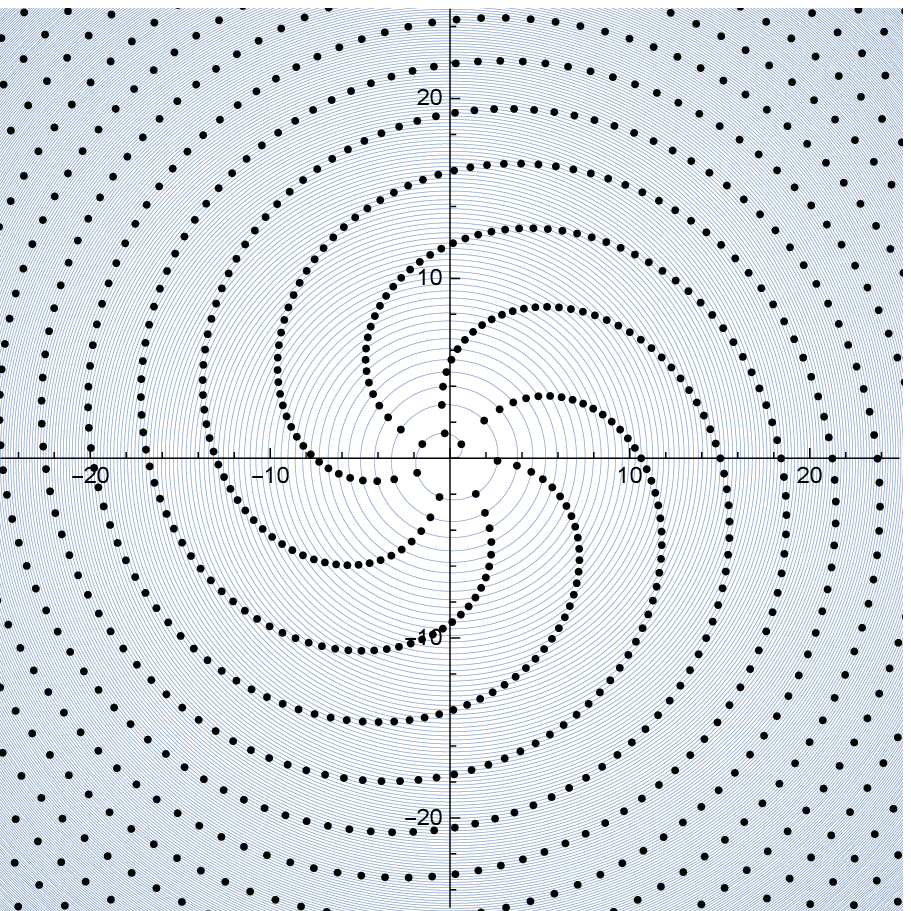}}
\end{subfigure}
\caption{Constant angular progressions on Fermat spiral\label{Fig}}
\end{figure}
In other words, we 
are interested in the sequence of points on the Fermat spiral 
that progresses by a constant angle $\alpha$ (see Figure \ref{Fig}).
A real number $\alpha$ is {\it badly approximable} if there exists a positive constant $C$
so that
$$
q|q \alpha - p|\ge C
$$
holds for all $(p,q)\in \Z\times \N$. It is known that $\alpha$ is badly approximable
if and only if partial quotients of continued fraction expansion of $\alpha$ is bounded by a positive integer 
$B$
(see \cite[Theorem 23]{Khinchin:64}). Indeed, $B$ and $C$ are roughly inverse-proportional;
\begin{equation}
\label{BC}
\frac{1}{B+2}\le C\le \frac 1{B},
\end{equation}
c.f. \cite[Theorem 1.9]{Bugeaud_Approx}. 
We use this notation $B$ and $C$ throughout this paper.
In particular, if $\alpha$ is a real quadratic irrational, then $\alpha$ is badly approximable,
due to Lagrange Theorem. In this note we will prove

\begin{thm}
\label{Main}
The following four statements are equivalent.
\begin{enumerate}
\item $X(\sqrt{t},\alpha)$ is relatively dense,
\item $X(\sqrt{t},\alpha)$ is uniformly discrete,
\item $X(\sqrt{t},\alpha)$ is a Delone set,
\item the angle $\alpha$ is badly approximable.
\end{enumerate}
\end{thm}

By definition of a Delone set, 
the proof is finished when we prove the equivalence between a) and d) in \S \ref{Dense}, 
and then the one between b) and d) in \S \ref{Discrete}. We prove them in a quantitative form.

\begin{thm}
\label{Qua}
Let $B, C$ be the above bounds of $\alpha$ when $\alpha$ is badly approximable.
\begin{enumerate}
\item If $\alpha$ is badly approximable, then $X(\sqrt{t},\alpha)$ is $6\sqrt{B}$-relatively dense.
\item If $X(\sqrt{t},\alpha)$ is $r$-relatively dense, then 
$\alpha$ is badly approximable with $B=\lfloor 96r^2 \rfloor$.
\item If $\alpha$ is badly approximable, then 
$X(\sqrt{t},\alpha)$ is $\sqrt{C}/2$-uniformly discrete. 
\item If $X(\sqrt{t},\alpha)$ is $s$-uniformly discrete, then $\alpha$ is badly approximable
with $C=s^2/53$.
\end{enumerate}
\end{thm}

Roughly speaking, the bound $B$ 
corresponds to $(r,s)$-Delone set where $r$ is proportional to $\sqrt{B}$ and $s$
is inverse-proportional to $\sqrt{B}$ as in (\ref{BC}).  
As $(r,s)$-Delone set is more uniform 
when $r/s$ is closer to $1$, 
the case $B=1$, i.e.,  $\alpha=(\sqrt{5}-1)/2$ is the best choice.
This shows an important role of the golden angle observed in many instances in phyllotaxis.

\section{Three Distance Theorem}
The torus $\T:=\R/\Z\simeq [0,1)$ is divided into subintervals by $$\{ n\alpha \pmod{1}\ |\
n\in \{0,1,\dots,N\}\}.$$ 
Steinhaus conjectured that there are at most three different lengths
appear in this partition for any $\alpha$ and $N$. This 
is first proved by S\'os \cite{Sos:57, Sos:58}
and then \'Swierczkowski
\cite{Swierczkowski:59}, Sur\'anyi \cite{Suranyi:58},
Halton \cite{Halton:65} and Slater \cite{Slater:67}.
We follow a formulation in Alessandri-Berth\'e \cite{Alessandri-Berthe:98}, see also
\cite[Chapter 2.6]{Allouche-Shallit:03} and \cite[Chapter 6.4, p.518]{Knuth:73}.
Let $\eta_{-1}=1$, $\eta_0=\alpha$ and 
$$
\eta_{k-1}=a_{k+1}\eta_k+\eta_{k+1}, \qquad 0\le \eta_{k+1}<\eta_k
$$
for $k\ge 0$ with a unique $a_{k+1}\in \N$.
When $\alpha$ is irrational, the sequence 
$(\eta_k)_{k=-1}^{\infty}$ 
is positive, strictly decreasing and given by Euclidean division process
of intervals. Concatenating
$$
\frac {\eta_{k-1}}{\eta_k}=a_{k+1}+\frac 1{\frac {\eta_{k}}{\eta_{k+1}}},
$$
the continued fraction expansion of 
$$
\alpha=[0;a_1,a_2,a_3,\dots]=\cfrac{1}{a_1+\cfrac{1}{a_2+\cfrac{1}{a_3+\cfrac{1}{\ddots}}}}
$$
is naturally associated with this algorithm. 
Define 
$p_{-1}=1,p_0=0,q_{-1}=0,q_0=1$ and 
$$
p_{k+1}=a_{k+1}p_k+p_{k-1},\qquad q_{k+1}=a_{k+1}q_k+q_{k-1}
$$
for $k\ge 0$. Then we observe that the convergent $[0,a_1,\dots,a_n]$ is equal to $p_n/q_n$
and $\eta_k=(-1)^{k}(q_k \alpha -p_k)$ by induction. 
We later use an important basic inequality 
\begin{equation}
\label{Below}
\eta_k> \frac{1}{q_{k+1}+q_{k}}>\frac 1{2q_{k+1}},
\end{equation}
see e.g. \cite[Chapter 2]{Rockett-Szusz:92} ,\cite[Theorem 13]{Khinchin:64},
\cite[Corollary 1.4]{Bugeaud_Approx}.
For a positive integer $N$ there exists a unique expression:
$$
N=s q_k+q_{k-1}+r
$$
with $1\le s\le a_{k+1}$ and $0\le r<q_k$. 
Three distance theorem reads 

\begin{thm}
\label{TDT}
The torus $\T$ is subdivided by $\{ n\alpha \pmod 1\ |\
n\in \{0,1,\dots,N\}\}$ into intervals of at most three lengths:
\begin{itemize}
\item $N+1-q_k$ intervals have length $\eta_k$,
\item $r+1$ intervals have length $\eta_{k-1}-s \eta_k$,
\item $q_k-(r+1)$ intervals have length $\eta_{k-1}-(s-1)\eta_k$.
\end{itemize}
\end{thm}

Rotating all points 
by $x\mapsto x+\ell \alpha \pmod{1}$ in $\T$ with a fixed $\ell\in\Z$, 
the same statement is valid for $\{n \alpha \pmod 1\ |\ 
n\in \{\ell,\ell+1,\dots, N+\ell\}\}$. 
See \cite{Sos:57, Sos:58, Swierczkowski:59, Suranyi:58, Halton:65, Slater:67} for the proof
of Theorem \ref{TDT}.
We just comment 
that it makes
easier to understand them after 
we comprehend the recursive (dynamical) induction structure of irrational rotation
$([0,1), x\mapsto x+\alpha)$. 
Inducing an irrational rotation $x\mapsto x+(-1)^k\eta_k$ on the interval $I$ of length $\eta_{k-1}$
to an interval of length $\eta_{k}$ located at the end of $I$ 
in the sign $(-1)^{k}$ direction, the first return map gives another irrational rotation 
$x\mapsto x+(-1)^{k+1}\eta_{k+1}$.
As the number of points $n\alpha \pmod{1}$ increases, we inductively observe the induced irrational rotations.
Three distance theorem can be easily understood as a consequence of this structure.

\begin{cor}
\label{Rat} 
$\alpha$ is badly approximable if and only if there is a positive constant $R$, that
for any $N\in \N$, the ratios of three distances generated by the $N+1$ points 
are bounded by $R$.
\end{cor}

\proof
Clearly $\eta_{k-1}-(s-1)\eta_k$ is the sum of the other two. We obtain
$$
1<\frac{\eta_{k-1}-(s-1)\eta_k}{\eta_k}\le \frac{\eta_{k-1}}{\eta_k} \le 1+a_{k+1}
$$
and
$$
1<\frac{\eta_{k-1}-(s-1)\eta_k}{\eta_{k-1}-s \eta_k}=1+\frac{\eta_k}{\eta_{k-1}-s \eta_k}
\le 1+\frac{\eta_k}{\eta_{k+1}} \le 2+a_{k+2}.
$$
On the other hand, 
for $n=a_{k}q_{k-1}+q_{k-2}+q_{k-1}-1=q_{k}+q_{k-1}-1$ 
there are exactly two lengths $\eta_{k-1}$ and $\eta_{k-2}-a_{k}\eta_{k-1}=\eta_{k}$
which appear in the partition.
From $\eta_{k-1}-a_{k+1}\eta_{k}=\eta_{k+1}$ we see
$$
|\eta_{k-1}/\eta_{k}-a_{k+1}|<1.
$$
Summing up, the bound $B$ of partial quotients
and $R$ satisfy
\begin{equation}
\label{BR}
B-1\le R\le B+2.
\end{equation}
\qed
%

\section{Diophantine condition for relative denseness}\label{Dense}
First we give a Diophantine characterization of relatively denseness. 
The distance of $x$ from $\Z$
is designated by
$\| x\|=\min_{i\in \Z}|x-i|$.
Hereafter we often use $$4\|x-y\| \le |\e(x)-\e(y)| \le 2\pi \|x-y\|$$
which follows from $|1-\e(x)|=2|\sin(\pi x)|$.

\begin{lem}
\label{Dio}
If $X(\sqrt{t},\alpha)$
 is $r$-relatively dense then
for any $m\ge 16r^2$ and $\beta\in [0,1)$, the system of inequalities
\begin{equation}
\label{Key}
m-3r\sqrt{m} \le n \le m+3r \sqrt{m}, \qquad \|n\alpha -\beta\|\le \frac{r}{2\sqrt{m}}.
\end{equation}
has a solution $n\in \N$. Conversely the system of inequalities (\ref{Key}) is solvable, 
then $X(\sqrt{t},\alpha)$ is $6r$-relatively dense.
\end{lem}

\proof
Assume that $X(\sqrt{t},\alpha)$ is $r$-relatively dense. Then 
for $z=\sqrt{m}\e(\beta)\in \C$ with $m\ge 16r^2$, we find 
$n\in \N$ that $\sqrt{n} \e(n \alpha)\in B(z,r)$.
Then
$$
r> |\sqrt{n}\e(n \alpha)-\sqrt{m}\e( \beta)|\ge |\sqrt{n}-\sqrt{m}|
$$
gives
$$
|n-m|\le (\sqrt{n}+\sqrt{m})r\le (2\sqrt{m}+r)r\le 3\sqrt{m}r.
$$
On the other hand,
\begin{eqnarray*}
r &\ge & |\sqrt{m}\e(n \alpha)-\sqrt{m}\e(\beta)|- |\sqrt{n}\e(n \alpha)-\sqrt{m}\e(n \alpha)|\\
&\ge& 4\sqrt{m}\|n\alpha -\beta\| -r
\end{eqnarray*}
shows
$$
\frac r{2\sqrt{m}} > \|n\alpha-\beta\|.
$$
Conversely, take a complex number $z=\sqrt{m}\e(\beta)$ with $m\ge 16r^2$ and
assume the system (\ref{Key}) has a solution $n\in \N$.
Then
\begin{eqnarray*}
&&|\sqrt{n} \e(n \alpha)-z|\\
&\le&
|\sqrt{n}-\sqrt{m}|+|\e(n \alpha)-\e(\beta)|\sqrt{m}\\
&\le&
\frac {|n-m|}{\sqrt{n}+\sqrt{m}} + 2\pi \sqrt{m} \| n\alpha-\beta\|
\le \frac{3\sqrt{m}r}{\sqrt{m}/2+\sqrt{m}} + 2\pi \sqrt{m} \frac{r}{2\sqrt{m}}< 6r.
\end{eqnarray*}
Here we used $n\ge m/4$ which follows from $m\ge 16r^2$.
This shows for any point $z\in \C$ with $m\ge 16r^2$, the intersection
$B(z,6r)\cap X(\sqrt{t},\alpha)$ is non empty.
For $m<16r^2$, we have $|\sqrt{m}\e(\beta)-\e(\alpha)|<4r+1$. By (\ref{rs}), 
$X(\sqrt{t},\alpha)$ is $6r$-relatively dense.
\qed
\bigskip



{\it Proof of a) and b) of Theorem \ref{Qua}.}
We show that if $\alpha$ is badly approximable with the constant $R$ in Corollary \ref{Rat}, then 
the system of inequalities (\ref{Key}) is solvable. 
Take $R$ in Corollary \ref{Rat} and fix $r>\max\{\sqrt{R/5},1\}$.  
For $m\ge 16r^2$, the closed interval
$[m-3r\sqrt{m}, m+3r\sqrt{m}]$ contains $N+1=\lfloor 3r\sqrt{m}\rfloor-\lceil -3r\sqrt{m}\rceil+1
$ integers.
By Theorem \ref{TDT}
the points $$n \alpha \pmod{1}\quad \text{ with }\quad n\in \Z \cap [m-3r\sqrt{m}, m+3r\sqrt{m}]$$
gives a partition of $\T$ into $N+1$ intervals of at most three lengths
where $$6r \sqrt{m}-1<N+1\le 6r\sqrt{m}+1.$$ 
Every subdivided interval has length not larger than $R/(6r \sqrt{m}-1)$,
because otherwise all lengths are larger than $1/(6r \sqrt{m}-1)$ which would 
cause overlapping.
Therefore any $\beta\in \T$ and there exists $n$ such that
$$
m-3r\sqrt{m}\le n \le m+3r\sqrt{m}, \qquad \| n \alpha -\beta \|\le \frac {R}{2(6r \sqrt{m}-1)}\le \frac r{2\sqrt{m}}
$$
which shows that (\ref{Key}) is solvable, recalling (\ref{rs}).
In view of Lemma \ref{Dio} and (\ref{BR}), to show
$r$-relatively denseness of $X(\sqrt{t},\alpha)$, we need
$$
r\ge 6 \max\left\{ \sqrt{\frac{B+2}5},1\right\}
$$
and therefore $r\ge 6\sqrt{B}$ suffices.

Let us assume that $X(\sqrt{t},\alpha)$ is $r$-relatively dense. By Lemma \ref{Dio}, there exists $r>0$
that for any $m\ge 16r^2$ and any $\beta\in \T$, there exists $n\in \N$ satisfying (\ref{Key}).
Our goal is to find an upper bound of $a_k$. 
We may assume $a_{k+1}\ge 2$ and select $m\ge 16r^2$ that
$$
[m-3r\sqrt{m}, m+3r \sqrt{m}] \cap \Z
$$
has cardinality $\lfloor a_{k+1}/2\rfloor q_k + q_{k-1}+t$ with $t\in \{1,2\}$
for 
\footnote{The inequality $q_{k+1}\ge 48r^2$ ensures the existence of $m$ not less than $16r^2$.} 
$q_{k+1}\ge 48r^2$.
This is possible since the function
$$
g(m)=\mathrm{Card}([m-3r\sqrt{m}, x+3r\sqrt{m}] \cap \Z)$$
from $\R_{\ge 0}$ to $\N\cup\{0\}$ satisfies
$$
g(0)=1,\quad \lim_{m\rightarrow \infty} g(m)=\infty\quad \text{ and }\quad |g(m+0)-g(m-0)|\le 2.
$$
By this choice of $m$, we have
\begin{equation}
\label{Width}
6r \sqrt{m}\ge \lfloor a_{k+1}/2\rfloor q_k+q_{k-1} \ge \frac {q_{k+1}}4.
\end{equation}
By Theorem \ref{TDT}, the torus $\T$ is partitioned into subintervals having 
three lengths $\eta_{k}$, $\eta_{k-1}-\lfloor a_{k+1}/2\rfloor \eta$
and $\eta_{k-1}-\lfloor a_{k+1}/2-1\rfloor \eta$ by $$
\{n\alpha \pmod{1} \ |\ m-3r\sqrt{m}\le n \le m+3r \sqrt{m}\}.$$
From the choice of $m$, Theorem \ref{TDT} and 
(\ref{Below}), a subinterval $I$ of 
length $$\eta_{k-1}-\lfloor a_{k+1}/2-1 \rfloor \eta_{k}\ge a_{k+1}\eta_k/2\ge a_{k+1}/(4q_{k+1})$$ 
exists.
Letting $\beta$ be the mid point of $I$, we have $\|n\alpha-\beta\|\ge a_{k+1}/(8q_{k+1})$.
Summing up, we obtain
$$
\frac{a_{k+1}}{8q_{k+1}}\le \frac r{2\sqrt{m}} \le \frac {12r^2}{q_{k+1}}
$$
by (\ref{Width}). This shows $a_{k+1}\le 96 r^2$ under the assumption 
$q_{k+1}\ge 48r^2$. Clearly $q_{k+1}<48r^2$ implies $a_{k+1}<48r^2$. 
Therefore we obtained $B\le 96r^2$.
\qed

\section{Uniform discreteness of $X(\sqrt{t},\alpha)$}\label{Discrete}

Finally we prove c) and d) of Theorem \ref{Qua}.
If $X(\sqrt{t},\alpha)$ is not $(s/2)$-uniformly discrete, then we can find $n,m\in \N$
with $n<m$ such that
$$
|\sqrt{n} \e(n\alpha)-\sqrt{m} \e(m\alpha)|<s.
$$
Similar computation with Lemma \ref{Dio} leads to
$$
|m-n|<2\sqrt{m}s, \qquad \|n\alpha-m\alpha\|\le \frac {s}{2\sqrt{m}}.
$$
Putting $k=m-n$, we have $
k\|k \alpha\|<s^2$. 
Therefore $\alpha$ is not badly approximable by the constant $C=s^2$. 

Conversely assume $\alpha$ is not badly approximable by the constant $C=s^2$, i.e., 
$\|k \alpha\| < \frac {s^2}{k}$ has an integer solution $k\in \N$. 
Then
\begin{eqnarray*}
&&|\sqrt{n} \e(n \alpha)-\sqrt{n+k}\e((n+k)\alpha)|\\
&\le& |\sqrt{n+k}-\sqrt{n}|+|\e(n \alpha)-\e( (n+k) \alpha)|\sqrt{n}\\
&\le& \frac k{\sqrt{n+k}+\sqrt{n}} + 2\pi \sqrt{n} \|k \alpha \|
\end{eqnarray*}

Putting $n=\lceil (k/s)^2 \rceil$ and using (\ref{rs}), the right side is estimated by
$$
s \frac k{\sqrt{k^2+s^2 k}+k} + 2\pi \sqrt{1+\frac {k^2}{s^2}} \frac {s^2}k \le \left(\frac 12+2\pi \sqrt{5}\right) s.
$$
This shows that $X(\sqrt{t},\alpha)$ is not $(1/4+\pi \sqrt{5})s$-uniformly discrete.

\section{Open Questions}

It is of interest to generalize this result to a strictly increasing function $f$ 
asymptotically equal to $\sqrt{t}$, or to 
substitute $\e(n\alpha)$ by a uniformly distributed sequence on $\T$. 

Is there a higher dimensional analogy of this result ? One possibility is to find a 
transitive action $T$ on $S^{d-1}$, so that 
$$
\left\{\left. \sqrt[d]{n}\ T^n(x_0)\ \right|\ n\in \N \right\}
$$
is a Delone set in $\R^d$.

\bigskip
\begin{center}
{\bf Acknowledgments}
\end{center}
\bigskip

The author came across this problem from an inspiring lecture delivered
by J.~F.~Sadoc at RIMS Kyoto in Oct. 2018. 
He wishes to show his heartfelt gratitude to
Y.~Yamagishi at Ryukoku Univ, who organized the
meeting, for advices, discussion 
and critical reading of the original manuscript. 
He is also indebted to J.~Lagarias, Y.~Bugeaud and A.~Haynes for crucial remarks and discussion. 
This research is partially supported by JSPS grants (17K05159, 17H02849, BBD30028). 

\bibliographystyle{amsplain}


\end{document}